\documentclass[12pt]{article}

\newtheorem{Le}{Lemma}
\newtheorem{Co}{Conjecture}
\begin{document}

\title{Correlation Inequality for Formal Series}
\date{February 10,  2013}
\author{Vladimir Blinovsky}
\date{\small
 Institute for Information Transmission Problems, \\
 B. Karetnyi 19, Moscow, Russia,\\
 Instituto de Matematica e Statistica, USP,\\
Rua do Matao 1010, 05508- 090, Sao Paulo, Brazil\\
vblinovs@yandex.ru}

\maketitle\bigskip

\begin{center}
{\bf Abstract}
\end{center}
We extend the considerations of the paper~\cite{2} and prove two correlation inequalities for totally ordered set. \bigskip

{\bf Introduction}
\bigskip

We extend the considerations of the paper~\cite{2} and prove two correlation inequalities (statement of Lemma below and inequality~(\ref{e2})) for totally ordered set. From other side in~\cite{1} was made two conjectures (statement of Lemma and inequality~(\ref{e2}))
for poset $2^X$ of subsets of finite set $X$ with FKG condition on probability measure (see~(\ref{e1})). In~\cite{2} was stated that considerations from it  lead to the proof of the Lemma under these FKG  conditions on measure $\mu$, but it turns out that  that considerations  are not sufficient
for the proof and the problem is still open. To solve these conjectures (as we show here it is sufficient to prove Lemma and then~(\ref{e2}) follows), if they are true, one need to make some additional efforts.

\bigskip

{\bf Main Text}
\bigskip

First we introduce class of correlation inequalities.

Assume that $f_1 ,\ldots ,f_n$ are nonnegative nondecreasing functions $2^X \to R.$ The expectation of a random variable $f: 2^X \to R$ with respect to $\mu$ we denote by $\langle f\rangle_\mu .$
For a subset $\delta\in [n]$ define
$$
E_\delta =\left\langle\prod_{i\in\delta} f_i \right\rangle_\mu .
$$
Let
$$
\sigma =\{ \sigma_1 ,\ldots ,\sigma_\ell \} 
$$
be a partition of $[n]$ into disjoint subsets. Define
$$
E_\sigma =\prod_{i=1}^{\ell} E_{\sigma_i} .
$$
Let $\lambda_1 =|\sigma_i | .$ We have $\sum_{i=1}^{\ell} \lambda_i =n.$ Let $\lambda (\sigma )=(\lambda_1 ,\ldots ,\lambda_{\ell})$ and $\lambda_1 \geq \ldots \geq\lambda_{\ell} .$
For a partition $\lambda$ of number $n$ define
$$
E_\lambda =\sum_{\sigma : \lambda (\sigma )=\lambda}E_{\sigma} .
$$
 We need  the following
\begin{Le}
\label{le2}
Consider the totally odered set $1,\ldots ,N$ with probability measure $\mu$ on it and let's functions $f_i ,\ i=1,\ldots , n$ are nonnegative and monotone nondecreasing.
Functional 
\begin{equation}
\label{e3}
E_n (f_1 ,\ldots ,f_n )=\sum_{\lambda\vdash n}c_{\lambda}E_{\lambda}
\end{equation}
where
$$
c_{\lambda}=(-1)^{\ell +1}\prod_{i=1}^{\ell}(\lambda_i -1)!
$$
is nonnegative.
\end{Le}
In~\cite{1} was conjectured that statement of Lemma (along with~(\ref{e2})) is true when probability measure $\mu$ on $2^X$ satisfies FKG conditions
\begin{equation}
\label{e1}
\mu (A\cap B)\mu (A\cup B)\geq \mu (A)\mu (B)
\end{equation}
and functions $f_i$ are nonnegative and monotone. 

 Note that under conditions from the Lemma in particular case $n=2$ Lemma gives Chebyshev inequality
  $$
  \langle f_1 f_2\rangle_\mu  \geq \langle f_1 \rangle_\mu \langle f_2 \rangle_\mu  
  $$
  Hence our proof can be considered as extention of Chebyshev inequality to multiple variables.
  For monotone functions $f_i (j),\ i=1,\ldots ,n;\ j=1,\ldots ,N$ we put 
\begin{equation}
\label{e100}
 f_i (1)=a_{i,1},\ f_i (j)=f_i (j-1)+a_{i,j},\ j=2,\ldots ,N,\ a_{i,j}\geq 0 .
 \end{equation} 
Then substituting in the formula 
\begin{equation}
\label{e21}
E_n (f_1 ,\ldots , f_n ) =\sum_{\lambda\vdash n}c_{\lambda} \sum_{\sigma :\ \lambda (\sigma )=\lambda}\prod_{i=1}^{\ell} \langle \prod_{j\in\sigma_i } f_j \rangle_{\mu}
\end{equation}
 coefficients $c_\lambda$ one can easily check that coefficient of the monomial 
$$
\prod_{j=0}^{N-1}\prod_{i=\sum_{s=1}^{j} m_s +1}^{\sum_{i=1}^{j+1}m_s} f_i (j+1)$$
in the rhs of~(\ref{e21}) is as follows
 \begin{equation}
 \label{e33}
 F_{m_1 ,\ldots ,m_N} (\mu )=\sum_{k_{i,j} :\ \sum_i ik_{i,j}=m_j}\prod_{j=1}^{N} \mu^{\kappa_j}(j)
 \frac{\prod_{j=1}^{N}m_j ! (-1)^{\sum_{j=1}^N \kappa_j -1}}{\prod_{j=1}^N \prod \left( i^{k_{i,j}}k_{i,j}!\right)},
 \end{equation}
 where $\kappa_j =\sum_i k_{i,j} .$ 
 Analogous monomials with $m_j$ factors $f_{i_p}(j)$ have the same coefficients.
 Indeed
 \begin{eqnarray*}
 && \sum_{\lambda\vdash n}c_{\lambda} \sum_{\sigma : \lambda (\sigma )=\lambda}\prod_{i=1}^{\ell}\langle \prod_{j\in\sigma_i}f_j \rangle \\
 &=& \sum_{\lambda\vdash n}\sum_{\sigma :\lambda (\sigma )=\lambda} (-1)^{\sum_{ij=1}^{N}\kappa_j -1}\prod ((i-1)! )^{\sum_{j=1}^N k_{i,j}}\prod_{i=1}^{\sum_{j=1}^{N}\kappa_j}\sum_{j=1}^{N}
 \mu (j) \prod_{s\in\sigma_i}f_s (j) ,
 \end{eqnarray*}
 where $\{ k_{i,j} \}$ are number of ocurence of the sets of cardinality $i$ in the projection 
 of partition $\sigma$ onto $m_j .$ Number of partitions of $\sigma$ with given $\{ k_{i,j} \}$ is
 $$
 \frac{\prod_{j=1}^N m_j !}{\prod_{j=1}^N \prod_i (i!^{k_{i,j}} k_{i,j}! )} .
 $$
 Using identity ([7, p.181])
 $$
 \sum_{ \{ k_i \} :\ \sum ik_i =n}\frac{y^k}{\prod (k_i ! i^{k_i})} =\frac{y (y+1)\ldots (y+n -1)}{n!}
 $$
 from~(\ref{e33}) we obtain formula
 \begin{eqnarray*}
 F_{m_1 ,\ldots ,m_N}(\mu )=&-& (-\mu (1))(1-\mu (1))\ldots (m_1 -1-\mu (1))\\
 &\times& (-\mu (2))(1-\mu (2))\ldots (m_2 -1-\mu (2))\ldots\\
 &\times& (-\mu (N))(1-\mu (N))\ldots (m_n -1-\mu (N))=-\prod_{j=1}^N \prod_{i=1}^{m_j} (i-1-\mu (j))  
 \end{eqnarray*}
 Using decomposition~(\ref{e100}) it is easy to see that to prove Lemma it is sufficient to prove the inequality 
 \begin{equation}
 \label{e200}
 \sum_{\{ m_s \}}F_{m_1 ,\ldots ,m_N}(\mu )\prod_{j=0}^{N-1}\prod_{i=\sum_{s=1}^j m_s +1}^{\sum_{s=1}^{j+1} m_s}\left(\sum_{t=1}^{j+1} a_{i,t}\right) \geq 0 .
 \end{equation}

 One can check that coefficient before the monomial
 $$
 \prod_{j=0}^{N-1} \prod_{i=\sum_{s=1}^j m_s +1}^{\sum_{s=1}^{j+1}m_s} a_{i,j+1}$$
  in the lhs of~(\ref{e200}) is
  \begin{equation}
  \label{e300}
  B(m_1 ,\ldots ,m_{N-1} )\stackrel{\Delta}{=} -\sum_{\{i_j \}}\prod_{j=1}^{N-1}{\sum_{s=1}^j m_s -\sum_{s=1}^{j-1} i_s \choose i_j} \prod_{i=1}^{i_j} (i-1-\mu (j))\times\prod_{i=1}^{n-\sum_{s=1}^{N-1}i_s }(i-1-\mu (N)) .
  \end{equation}
  Thus to prove~(\ref{e200}) and complete the proof of Lemma it is sufficient to prove the inequality
  \begin{equation}
  \label{e300}
  B(m_1 ,\ldots ,m_{N-1})\geq 0.
  \end{equation}
  We prove this inequality by induction on $m_j$.
 Let's~(\ref{e300}) is true  for $m_N =n-\sum_{s=1}^{N-1}m_s \geq 1,$ then for $m_N +1$ we have the expression for $B(m_1 ,\ldots ,m_{N-1})$
 (we use the identity ${\ell\choose p}={\ell -1\choose p}+{\ell -1\choose p-1} $)
 \begin{eqnarray}
 \label{e400}
 &-&\sum_{\{ i_s \}}\prod_{j=1}^{N-1} {\sum_{s=1}^{j}m_s -\sum_{s=1}^{j-1}i_s \choose i_j}
 \prod_{j=1}^{i_j} (i-1-\mu (j))\\
 &\times& \prod_{i=1}^{n-\sum_{j=1}^{N-1}i_j}(i-1-\mu (N))(n-\sum_{s=1}^{N-1}i_s -\mu (N))  .
 \end{eqnarray}
 Last expression is nonnegative due to the induction proposal and the fact that (because $m_N \geq 1$)
 $$
 n-\sum_{s=1}^{N-1}i_s -\mu_N \geq 0 .
 $$
 Step by step using induction we come to the expression for $B$ with $m_N =1$ (we assume at first that $m_N >0$) and start induction on $m_{N-1}.$ Let~(\ref{e300}) is true for $m_N$, then the expression~(\ref{e400}) for $m_N +1$ 
 \begin{eqnarray}\label{er}
 &&-\sum_{\{ i_s \}}\prod_{j=1}^{N-1}{\sum_{s=1}^{j}m_s -\sum_{s=1}^{j-1}i_s \choose i_s}\prod_{j=1}^{i_s}(i-1-\mu (j))\\
 &\times& \prod_{i=1}^{n-\sum_{s=1}^{N-1}i_s}(i-1-\mu (N))\left( n-\sum_{s=1}^{N-2} i_s -\mu (N-1)-\mu (N)\right) .\nonumber
 \end{eqnarray}
 Because $m_N =1$, then
 $$
 n-\sum_{s=1}^{N-2}i_s -\mu (N-1)-\mu (N) \geq 0
 $$
 and thus  by induction hypothesis  expression~(\ref{er}) is nonnegative.
 Continuing this process to other $m_j ,\  j=N-3,\ldots ,1$ we come to the situation when $m_{N_0} =1,
 m_j =0,  j\neq N_0$ for some $N_0 \in [N].$ Thus to complete the induction we need to prove that expression~(\ref{e200}) in the case, when $m_j =1$ only for one value of $j,$ and all other $m_s =0$. But  this negativeness immideatly follows  from the relation
 $$
 B(0,\ldots , 0,1,0,\ldots ,0)=\mu (j) .
 $$
 This proves Lemma.   
  \bigskip
  
 Next we consider the set of formal series $P[[t]]$, whose coefficients are monotone nondecreasing nonnegative functions on $2^X.$ Then $p(A)=p_1 (A)t +p_2 (A)t^2 +\ldots \in P[[t]]$. In~\cite{1} was formulated the following
\begin{Co}
\label{co2}For FKG probability measure $\mu$
the following inequality is true
\begin{equation}
\label{e2}
1-\prod_{A\in 2^X}(1-p(A))^{\mu (A)}\geq 0.
\end{equation}
\end{Co}
The inequality~(\ref{e2}) is understood as non negativeness of coefficients of formal series obtained by series expansion of the product on the left- hand side of this inequality. 

We will prove, that inequality~(\ref{e2}) follows from inequalities
\begin{equation}
\label{e56}
E_n (f_1 ,\ldots, f_n)\geq 0 
\end{equation}
for all $n$ and hence it is sufficient to prove last inequalities and then inequality~(\ref{e2}) follows under the same conditions on $\mu$.

We make some transformations of the expression in the lhs of~(\ref{e2}).
We have
\begin{eqnarray*}
&&
1-\prod_{A\in 2^X}(1-p(A))^{\mu (A)}=1-\exp\left\{\left\langle \ln (1-p)\right\rangle_\mu\right\}\\
&=&1-\exp\left\{ -\sum_{i=1}^{\infty}\frac{1}{i}\langle p^i \rangle_{\mu}\right\}\\
&=& \sum_{j=1}^{\infty}\frac{(-1)^{j+1}}{j!}\left(\sum_{i=1}^{\infty}\frac{1}{i} \langle p^i \rangle_{\mu}\right)^j \\
&&= \sum_{j=1}^{\infty}\frac{(-1)^{j+1}}{j!}\sum_{\{ q_s\} :\sum q_s =j}{j\choose q_1 ,\ldots , q_j}\sum_{\{ i_s \}}
\frac {\prod_{s=1}^j \langle p^{i_s} \rangle^{q_s}_\mu}{(i_1)^{q_1} (i_2)^{q_2}\ldots (i_j)^{q_j}} .\end{eqnarray*}
Next remind that the number of partitions of $n$ with given set  $\{ q_i\}$ of occurrence of $i$ is equal to
$$
\frac{n!}{\prod_i (i!)^{q_i}q_i !} .
$$
Continuing the last chain of identities and using last formula we obtain
\begin{eqnarray}
\label{e4}
 &&1-\prod_{A\in 2^X}(1-p(A))^{\mu (A)}\\ \nonumber
  &&=\sum_{n=1}^{\infty}\frac{1}{n!} \sum_{\lambda\vdash n}\sum_{\sigma :\ \lambda (\sigma )=\lambda} (-1)^{\sum_i q_i +1}\frac{n! \prod ((i-1)!)^{q_i}}{\prod_i (i!)^{q_i} q_i !} \prod_i \langle p^{i} \rangle^{q_i}_\mu \\
  && \nonumber
  =\sum_{n=1}^{\infty}\frac{1}{n!}   \sum_{\lambda\vdash n}  (-1)^{\ell (\lambda )+1}\prod_i (\lambda_i -1)!
  \sum_{\sigma :\ \lambda (\sigma )=\lambda} 
 E_\sigma (p,\ldots ,p)\\
 && \nonumber
 =  \sum_{n=1}^{\infty}\frac{1}{n!}   \sum_{\lambda\vdash n} c_\lambda E_\lambda (p, \ldots ,p)\\
 && \nonumber
 =\sum_{n=1}^{\infty}\frac{1}{n!}E _n (p,\ldots ,p )
  \end{eqnarray}
  Hence now to prove the conjecture~\ref{co2} we need to show that
  \begin{equation}
  \label{et}
  E_n (p,\ldots ,p )\geq 0.
  \end{equation}
  But  the coefficients of the formal series 
  $E_n (p,\ldots p)$ are the sums of $E_n (p_{i_1} ,\ldots , p_{i_n})$ for  multisets
  $\{ i_1 ,\ldots ,i_n \} .$ This completes the proof that inequality~(\ref{e2}) follows from inequalities~(\ref{e56}) under the same conditions on $\mu .$
  
  Thus because we prove Lemma, we prove inequality~(\ref{e2}) for totally ordered lattice and this is our main result.
  
 \bigskip
  
{\bf Remark}
 \bigskip
 
To extend Lemma for the conditions~(\ref{e1})  one can try to find proper expansion for the monotone functions $f_i$ which extend expansion~(\ref{e100}) to the case of poset $2^X .$
   

\begin{thebibliography}{99}
\bibitem{2} Blinovsky, V.,  A Proof of One Correlation Inequality, Problems of Inform. Transm., 2009, vol. 45, no. 3, pp. 264- 269.
\bibitem{1}  Sahi, S., Higher Correlation Inequalities, Combinatorica, 2008, vol.28, no. 2, pp. 209- 227.
\bibitem{3} Riordan J., Combinatorial identities, Wiley, Ney York, 1968
\end{thebibliography}
\end{document}